\documentclass[letterpaper,10pt,reqno]{amsart}
\usepackage{indentfirst} % Indent first paragraph after section header
\usepackage{amssymb}
\usepackage{mathtools} % improves amsmath
\usepackage{mathabx} % makes many symbols (as $\leq$) more beautiful
\usepackage{amsthm}
\usepackage{thmtools}
\usepackage{enumitem} % special itemize with custom tags and labels that work
\usepackage[backref=page,colorlinks=true]{hyperref} % Links in the pdf. Backref option: each bibliographical entry denotes where it was cited
\usepackage[usenames,dvipsnames]{xcolor} % With XCOLOR, printed quality is good, (but with COLOR it's bad)
\usepackage{ifthen}
\usepackage{tikz}
\usetikzlibrary{decorations.pathreplacing}
\usepackage{tikz-cd}
  
\renewcommand*{\backref}[1]{}
\renewcommand*{\backrefalt}[4]{\tiny
  \ifcase #1 (\textbf{NOT CITED.})%
  \or    (Cited on page~#2.)%
  \else   (Cited on pages~#2.)%
  \fi}

\makeatletter

\def\MRbibitem{\@ifnextchar[\my@lbibitem\my@bibitem}

\def\mybiblabel#1#2{\@biblabel{{\hyperref{http://www.ams.org/mathscinet-getitem?mr=#1}{}{}{#2}}}}

\def\myhyperanchor#1{\Hy@raisedlink{\hyper@anchorstart{cite.#1}\hyper@anchorend}}

\def\my@lbibitem[#1]#2#3#4\par{%
  \item[\mybiblabel{#2}{#1}\myhyperanchor{#3}\hfill]#4%
  \@ifundefined{ifbackrefparscan}{}{\BR@backref{#3}}%
  \if@filesw{\let\protect\noexpand\immediate% write to aux-file
    \write\@auxout{\string\bibcite{#3}{#1}}}\fi\ignorespaces%
}

\def\my@bibitem#1#2#3\par{%
  \refstepcounter\@listctr% standard tex item counter for the generic item number
  \item[\mybiblabel{#1}{\the\value\@listctr}\myhyperanchor{#2}\hfill]#3%
  \@ifundefined{ifbackrefparscan}{}{\BR@backref{#2}}%
  \if@filesw\immediate\write\@auxout% write to aux-file
    {\string\bibcite{#2}{\the\value\@listctr}}\fi\ignorespaces%
}

\makeatother
\DeclareFontFamily{U} {MnSymbolA}{}
\DeclareFontShape{U}{MnSymbolA}{m}{n}{
   <-6> MnSymbolA5
   <6-7> MnSymbolA6
   <7-8> MnSymbolA7
   <8-9> MnSymbolA8
   <9-10> MnSymbolA9
   <10-12> MnSymbolA10
   <12-> MnSymbolA12}{}
\DeclareFontShape{U}{MnSymbolA}{b}{n}{
   <-6> MnSymbolA-Bold5
   <6-7> MnSymbolA-Bold6
   <7-8> MnSymbolA-Bold7
   <8-9> MnSymbolA-Bold8
   <9-10> MnSymbolA-Bold9
   <10-12> MnSymbolA-Bold10
   <12-> MnSymbolA-Bold12}{}
\DeclareSymbolFont{MnSyA} {U} {MnSymbolA}{m}{n}
 \DeclareFontFamily{U} {MnSymbolC}{}
\DeclareFontShape{U}{MnSymbolC}{m}{n}{
  <-6> MnSymbolC5
  <6-7> MnSymbolC6
  <7-8> MnSymbolC7
  <8-9> MnSymbolC8
  <9-10> MnSymbolC9
  <10-12> MnSymbolC10
  <12-> MnSymbolC12}{}
\DeclareFontShape{U}{MnSymbolC}{b}{n}{
  <-6> MnSymbolC-Bold5
  <6-7> MnSymbolC-Bold6
  <7-8> MnSymbolC-Bold7
  <8-9> MnSymbolC-Bold8
  <9-10> MnSymbolC-Bold9
  <10-12> MnSymbolC-Bold10
  <12-> MnSymbolC-Bold12}{}
\DeclareSymbolFont{MnSyC} {U} {MnSymbolC}{m}{n}

\DeclareMathSymbol{\top}{\mathord}{MnSyA}{219} % smaller symbol for transpose
\DeclareMathSymbol{\plus}{\mathord}{MnSyC}{20} % a smaller plus sign

\declaretheorem[numberwithin=section]{theorem}
\declaretheorem[sibling=theorem]{lemma}

\hypersetup{bookmarksdepth = 3} % Depth of sections/subs... to have bookmark links in the pdf;
\numberwithin{equation}{section}     % Makes labeled equations easier to find.

\setlist[enumerate,1]{label={\upshape(\alph*)},ref=\alph*}
\setlist[enumerate,2]{label={\upshape(\arabic*)},ref=\arabic*}

\def\phi{\varphi}

\newcommand{\vertiii}[1]{{\left\vert\kern-0.25ex\left\vert\kern-0.25ex\left\vert #1 
    \right\vert\kern-0.25ex\right\vert\kern-0.25ex\right\vert}}
\newcommand{\invertiii}[1]{{\vert\kern-0.25ex\vert\kern-0.25ex\vert #1 
    \vert\kern-0.25ex\vert\kern-0.25ex\vert}}

\begin{document}

\title{On Lie Groups and The Theory of Complex Variables}
\date{\today}

%\subjclass[2010]{37D35, 37A10, 37A35}

\begin{thanks}
{I am much indebted to Godofredo Iommi Echeverría for a careful  reading of a first
version  of this note.}
\end{thanks}

\author[G.~Iommi Amun\'ategui]{Godofredo Iommi Amun\'ategui}
\address{Instituto de F\'isica, Pontificia Universidad Cat\'olica de Valpara\'iso (P.UCV)}
\email{\href{godofredo.iommi@pucv.cl}{godofredo.iommi@pucv.cl}}
%\urladdr{\url{http://http://www.mat.uc.cl/~giommi/}}

\begin{abstract}
In this note we envisage the relation existing between the Lie Groups and the Theory of Complex Variables. In particular, it is shown that the dimensions of the irreducibles representations of $SU(N)$ may be written in terms of the Eisenstein integers and an identity is built up between the  imaginary parts of the dimensions of the irreducible representations of the Lie Groups $SU(3)$ and $Sp(4)$. 
\end{abstract}

\maketitle

\section{Introduction}
In a work dealing with the classes of binary quadratic forms with complex integral coefficients, G. Eisenstein \cite{d} introduced the numbers $a+b \omega$, where $a$ and $b$ are real integers and $\omega$ is an imaginary cube roof of unity.

D. Speiser \cite[Fig.7 and Fig. 22]{s}  noticed a curious connection between Lie Groups of rank $2$ and the theory of complex variables. In particular, he pointed out that in the lattice formed by the dimension function of $SU(3)$, the values are arithmetical means of their closed neighbors, fact which is reminiscent of harmonic functions. Consequently, he proposed the following identity for the dimension $\text{Dim}(p_1, p_2)$ of an irreducible representation $(p_1, p_2)$ of $SU(3)$\footnote{To obtain equation \eqref{spe} use must be made of 
$z= \left(p_1 + \frac{p_2}{2} \right) +i p_2 \frac{\sqrt{3}}{2}.$}:
\begin{equation}  \label{spe}
\text{Dim}(p_1, p_2)= \frac{1}{2} p_1 p_2 (p_1 + p_2) = \frac{1}{3 \sqrt{3}} \text{Im}z^3.
\end{equation}

In Section \ref{sec:2} we shall clarify the meaning of this equation and  subsequently we shall write the irreducible representation of $SU(N)$ in terms of the Eisenstein  numbers. In Section \ref{sec:3} an identity is built up between the imaginary parts of the dimensions of the irreducible representations of the Lie Groups $SU(3)$ and $Sp(4)$.

\section{Eisenstein integers and the dimensions of the irreducible representations of the unitary groups} \label{sec:2}

The Eisenstein integers numbers are $a+b \omega$, where $\omega= \frac{-1 + i \sqrt{3}}{2}$ is one of the cubic roots of unity, and the others are $1$ and $ \omega^2=\frac{-1 - i \sqrt{3}}{2}$. These numbers form a triangular lattice \cite[Fig 8.10 (a), p. 221)]{cg} and \cite[pp. 2--5]{c}. 

\begin{lemma}
The dimensions of the irreducible representations $(p_1, p_2)$ of $SU(3)$ may be written through Eisenstein numbers as:
\begin{equation} \label{2.1}
\text{Dim}(p_1,p_2)= \frac{\text{Im}(a+ b \omega )^3}{3 \sqrt{3}}.
\end{equation}
\end{lemma}

The imaginary part of $(a+b\omega)^3$ is:
\begin{equation} \label{2.3}
\text{Im}(a+ b \omega )^3=3ab(a\omega + b \omega^2) = i 3 \sqrt{3} \frac{(a-b)ab}{2}.
\end{equation}

The number $N(a,b)= \frac{1}{3 \sqrt{3}} \text{Im}(a+b \omega)^3 =\frac{(a-b)ab}{2}$ may be assigned to each Eisenstein lattice point. If we set $a=p_1 +p_2$ and $b=p_2$, we recover $D(p_1, p_2)$.

Long ago H. Weyl \cite{w} obtained the branching law for the groups of linear transformations. Hereafter we shall restrict his result to the case of unitary groups. For briefness we omit the details and we state the reduction $SU(N)\rightarrow SU(N-1)$: the $SU(N)$ irreducible representation $(p_1, p_2, \dots, p_{N-1})$ reduces into $SU(N-1)$ irreducible representations according to the formula

\begin{align} \left[ \left(	p_1, p_2, \dots, p_{N-1}\right)	\right] = \nonumber
\sum_{k_{1,1}=1}^{p_1} \sum_{k_{1,2}=1}^{p_2} \cdots  &\\ \sum_{k_{1,N-1}=1}^{p_{N-1}} 
\left(p_{1}-k_{1,1}+ k_{1,2}, p_2 -k_{1,2}+ k_{1,3}, \dots, p_{N-2} - k_{1, N-2} + k_{1, N-1}		\right).\label{2.4}
\end{align}

In order to write the dimension of $(p_1, p_2, \dots, p_{N-1})$ in terms of Eisenstein numbers the reduction chain $SU(N) \rightarrow SU(N-1) \rightarrow \dots \rightarrow SU(3)$ must be considered, i.e. the $SU(3)$ content of  $(p_1, p_2, \dots, p_{N-1})$ must be displayed. As an example, let us examine the case of an irreducible representation $SU(5)$. Taking into account the chain $SU(5) \rightarrow SU(4) \rightarrow  SU(3)$, the dimension of $(p_1, p_2, \dots, p_{N-1})$ may be expressed as follows:

\begin{eqnarray*} \label{2.5}
\text{Dim}(p_1, p_2, p_3, p_4)= 
\sum_{k_{1,1}=1}^{p_1} \sum_{k_{1,2}=1}^{p_2} \sum_{k_{1,3}=1}^{p_3} \sum_{k_{1,4}=1}^{p_4} \sum_{k_{2,1}=1}^{p_{1}-k_{1,1}+k_{1,2}} \sum_{k_{2,2}=1}^{p_{2} -k_{1,2}+ k_{1,3}} \sum_{k_{2,3}=1}^{p_{3}-k_{1,3}+k_{1,4}}  &\\
\frac{1}{3 \sqrt{3}} I_m\Big((a-b-k_{1,1}+ k_{1,2}-k_{2,1}+ k_{2,2}) + &\\ \frac{1}{2} (b-k_{1,2}+k_{1,3}-k_{2,2}+k_{2,3}) +i \frac{\sqrt{3}}{2} (b-k_{1,2}+k_{1,3}-k_{2,2}+ k_{2,3}) 	\Big)^3
\end{eqnarray*}

Although such identities acquire an involved aspect their underlying structure is transparent. The general result is:

\begin{lemma} \label{lem:2}
The dimension of the irreducible representation  $(p_1, p_2, \dots, p_{N-1})$ of $SU(N)$ may be written in terms if the Eisenstein numbers as:
\begin{eqnarray*} \label{2.6}
\text{Dim}(p_1, p_2, p_3, p_4)= &\\
\sum_{k_{1,1}=1}^{p_1} \cdots \sum_{k_{1, N-1}=1}^{p_{N-1}} \sum_{k_{2,1}=1}^{p_1-k_{1,1}+k_{1,2}} \cdots &\\ \sum_{k_{2,N-2}=1}^{p_{N-2}-k_{1, N-2}+k_{1,N-1}} \sum_{k_{3,1}=1}^{p_{1} -k_{1,1}+ k_{1,2}-k_{2,1}+k_{2,2}} \cdots \sum_{k_{3, N-3}=1}^{p_{1}-k_{1,N-3}+k_{1,N-2}-k_{2, N-3}+ k_{2, N-1}} \cdots   &\\
\frac{1}{3 \sqrt{3}} I_m\Big((a-b-k_{1,1}+ k_{1,2}-k_{2,1}+ k_{2,2} -k_{3,1}+ k_{3,2}- \dots) + &\\ \frac{1}{2} (b-k_{1,2}+k_{1,3}-k_{2,2}+k_{2,3}-k_{3,2}+k_{3,3}-\dots ) &\\ +i \frac{\sqrt{3}}{2} (b-k_{1,2}+k_{1,3}-k_{2,2}+ k_{2,3}- k_{3,2}+k_{3,3}- \dots) 	\Big)^3.
\end{eqnarray*}
\end{lemma}

Remark that, for $N>3$, this formula consists of $\frac{1}{2} (N+2)(N-3)$ summations.

\section{Concerning an identity on the complex plane}\label{sec:3}

%Hereafter we pursue an inquiry  about the relationship existing between  
%the Lie Groups of rank $2$ and the Theory of the Complex Variables Previously
% we have pointed out an equality between  the dimensions of the irreducible
%representations of the Unitary Group $SU(3)$, assigned to the lattice point $( p_1 , p_2  )$
%and the imaginary part of $z^3$ with
%
%\begin{equation} \label{1}
%z= \left(p_1 + \frac{p_2}{2} \right) + i \left(p_2   \frac{\sqrt{3}}{2} \right).
%\end{equation}
%It turns out that
%\begin{equation} \label{2}
%\text{Dim} (p_1, p_2) = \frac{1}{2!} p_1 p_2 \left(p_1 + p_2 \right) = \frac{1}{3 \sqrt{3}} \text{Im} z^3.
%\end{equation}
%(For details see \cite{1}).

In this section our purpose is to build up an identity between the imaginary parts of the  dimensions of the irreducible representations of the  Unitary  Group $SU(3)$ and the dimensions of the irreducible representations of the Symplectic Group $Sp(4)$. These  Groups are subgroups of the Unitary  Group $SU(4)$. The Group $SU(4)$ has rank $3$, hence its lattice is $3$-dimensional and to each lattice point $(p_1, p_2 , p_3  )$ corresponds an  irreducible representation  whose dimension is given by :
\begin{equation} \label{3}
\text{Dim} ( p_1,p_2 ,p_3 )= \frac{1}{2! 3!} p_1 p_2 p_3 (p_1 + p_2)(p_2 + p_3)(p_1 +p_2 + p_3 ) \end{equation}
where $p_1, p_2$  and  $p_3$    are positive integers.

In order to construct   such an identity ,we shall follow a procedure whose main steps are:

\begin{enumerate}
\item[(a)]  \label{a}The  branching rule for the reduction  $SU(4)  \rightarrow SU(3)$.

\item[(b)]  The branching rule for the reduction $SU(4)   \rightarrow Sp(4) $   and the geometrical transformation which allows us to express the $Sp (4)$ lattice point $( q_1, q_2)$ by means of the complex variable $z'$.\label{b}
\end{enumerate}
The final move is nothing but an adequate combination
of $(a)$ and $(b)$.

%Long ago Hermann Weyl established a branching rule for the Group of unimodular
%linear transformations in $N$ dimensions \cite{2}. His result is valid, as well, for the
%unitary case. Dressed, so to say, in a new notation, it may be written for the case
%$SU(4)  \rightarrow      SU(3)$  as:

For $N=4$, Weyl's branching rule \ref{2.4} may be written as:
\begin{equation} \label{3.2}
[(p_1, p_2, p_3)] = \sum_{a_1=1}^{p_1}\sum_{a_2=1}^{p_2}\sum_{a_3=1}^{p_3} (p_1 -a_1 +a_2, p_2-a_2+a_3),
\end{equation}
where $( p_1 , p_2   , p_3  )$ corresponds to an irreducible representation of $SU(4)$  in the Group lattice. The square brackets indicate the $SU(3)$ content of $(p_1 ,p_2  ,p_3  )$.

At this point it seems in order to recall that the Symplectic Group $Sp (4)$ consists of
the subset of unitary matrices in $4$ dimensions $(4=2l)$ which leaves invariant a skew
symmetric bilinear form:
\begin{equation*}
\sum_{i=1}^2 (x_i y_{i+2}-x_{i+2}y_i).
\end{equation*}

The existence of a non-degenerate skew symmetric form requires an even number of
dimensions. Besides, let us remark that the reduction  $SU(4)  \rightarrow     Sp(4)$  has been solved using a geometrical method \cite{i}. This reduction admits three cases which take into account the $Sp(4)$ lattice diagram symmetry:

\begin{equation} \label{3.3}
p_1 < p_3,   \quad [(p_1, p_2, p_3)] = \sum_{\rho=0}^{p_1-1}\sum_{\lambda=0}^{p_2-1} (p_1 + p_3-1-2p, 1+\rho+ \lambda) 
\end{equation}

\begin{equation} \label{3.4}
p_1 = p_3,   \quad [(p_1, p_2, p_3)] = \sum_{\rho=0}^{p_1-1}\sum_{\lambda=0}^{p_2-1} (2 p_1 -1-2 \rho, 1 + \rho + \lambda) 
\end{equation}
(or the same expression with $p_3$ instead of $p_1$).
\begin{equation} \label{3.5}
p_1 > p_3,   \quad [(p_1, p_2, p_3)] = \sum_{\rho=0}^{p_3-1}\sum_{\lambda=0}^{p_2-1} (p_1 + p_3-1-\rho, 1 + \lambda + \rho).
\end{equation}

A geometric consideration of the  $2-$dimensional lattice of $Sp(4)$ is crucial  to  find
a  transformation    which  gives  room  to an identity  between  the  irreducible 
representations of $Sp(4)$ and the imaginary part of a complex expression.  To   achieve
our goal  let us envisage the $Sp (4)$ lattice point $(q_1 ,q_2 )$ in such a manner that the   $q_1-$axis coincides with the $x-$axis of the complex plane. We get:

\begin{equation} \label{3.6} 
 z'= -\frac{q_1}{2} +i \left(\frac{q_1}{2}	+ q_2	\right). 
\end{equation}
From this  transformation and the dimension formula for  the  irreducible representations  of $Sp (4)$
\begin{equation} \label{3.7} 
\text{Dim} (q_1 ,q_2 ) = \frac{1}{6} q_1 q_2 (q_1 + q_2)(q_1 + 2 q_2).
\end{equation}
we deduce that
\begin{equation} \label{3.8} 
\text{Dim} ( q_1 , q_2) = \text{Im} \frac{z'^4}{6}.
\end{equation}

Let    $ [(p_1, p_2, p_3)]_{SU(3)}$ and    $ [(p_1, p_2, p_3)]_{Sp(4)}$            denote     respectively \eqref{3.2}, \eqref{3.3},\eqref{3.4} and \eqref{3.5}. We may state  the identity on the complex plane in the symbolic form:

\begin{lemma}\label{I}
\begin{equation*}
\text{Im} \frac{z^3}{3 \sqrt{3}} [(p_1, p_2, p_3)]_{SU(3)} = \text{Im} \frac{z'^4}{6} [(p_1, p_2, p_3)]_{Sp(4)}.
\end{equation*}
where  to  each resulting  term of the  decompositions  must be applied  the 
corresponding coordinate transformation.
\end{lemma}

To illustrate Lemma \ref{I}, let us work out the reduction  corresponding to
the irreducible representation $(5,3,2)$ of dimension $1000$ of $SU(4)$ :
\begin{align*}
\text{Im} \frac{z^3}{3 \sqrt{3}}  \left( \sum_{a_1=1}^{5}\sum_{a_2=1}^{3}\sum_{a_3=1}^{2} (5 -a_1 +a_2, 3-a_2+a_3)		\right)^3 = &\\ \text{Im} \frac{z'^4}{6} \left( \sum_{\rho=0}^{1}\sum_{\lambda=0}^{2} (6-2\rho, 1 + \lambda + \rho) \right)^4.
\end{align*}

Through  algebraic	manipulations, finally,  we obtain:

\begin{align*} \label{9}
 \frac{1}{3 \sqrt{3}} \text{Im} \Big( \Big[\left(\frac{5}{2} + i3 \frac{\sqrt{3}}{2} \right)^3	
 + 2 	\left(\frac{7}{2} + i3 \frac{\sqrt{3}}{2} \right)^3 + 2\left(\frac{9}{2} + i3 \frac{\sqrt{3}}{2} \right)^3 &\\
+2  \left(\frac{11}{2} + i3 \frac{\sqrt{3}}{2} \right)^3 + 2\left(\frac{13}{2} + i3 \frac{\sqrt{3}}{2} \right)^3 + \left(\frac{15}{2} + i3 \frac{\sqrt{3}}{2} \right)^3  \Big]  &\\
+\Big[\left(3 + i 2 \sqrt{3} \right)^3 + 	\left(4 + i 2 \sqrt{3} \right)^3	+ \left(5 + i 2 \sqrt{3} \right)^3 + \left(6 + i 2 \sqrt{3} \right)^3 + \left(7 + i 2 \sqrt{3} \right)^3	\Big]+ &\\ \Big[   \left(3 + i \sqrt{3} \right)^3    +  2\left(4 + i \sqrt{3} \right)^3 +  2\left(5 + i \sqrt{3} \right)^3
+ 2\left(6 + i\sqrt{3} \right)^3+  2\left(7 + i \sqrt{3} \right)^3 + &\\ \left(8 + i \sqrt{3} \right)^3\Big]+
\Big[\left(\frac{7}{2} + i \frac{\sqrt{3}}{2} \right)^3 +\left(\frac{9}{2} + i \frac{\sqrt{3}}{2} \right)^3 
+\left(\frac{11}{2} + i \frac{\sqrt{3}}{2} \right)^3+\left(\frac{13}{2} + i \frac{\sqrt{3}}{2} \right)^3 &\\
+\left(\frac{15}{2} + i \frac{\sqrt{3}}{2} \right)^3\Big] \Big) =&\\ 
\frac{\text{Im}}{6} \Big(\Big[(4i-3)^4 + (5i-3)^4 +(6i-3)^4 +(4i-2)^4 +(5i-2)^4+(6i-2)^4 \Big] \Big).
\end{align*}

Dimensional verification:
\begin{align*}
\text{Dim} [(5,3,2)]_{SU(3)}=[6+30+54+84+120+81] + [10+24+42+64+90] +&\\ [8+30+48+70+96+63]+ [6+10+15+21+28].
\end{align*}

\begin{align*}
\text{Dim} [(5,3,2)]_{Sp(4)}=[56+160+324+64+140+256].		
\end{align*}

\end{document}